\documentclass[a4paper,11pt]{amsart}
\usepackage{color}
\usepackage{amssymb}
\usepackage{mhequ}
\usepackage{cite}

\def \Reel{\mathbb{R}}

\def \Nat{\mathbb{N}}

\newcommand{\dte}%{o(\delta^2)??}%
{\delta^{4}}
\newcommand{\dnte}%{o(\delta^2_n)??}%
{\delta_n^{4}}

\DeclareFontFamily{OT1}{rsfs}{}
\DeclareFontShape{OT1}{rsfs}{m}{n}{ <-7> rsfs5 <7-10> rsfs7 <10->
rsfs10}{} \DeclareMathAlphabet{\mathscr}{OT1}{rsfs}{m}{n}
\newcommand{\bel}[1]{\begin{equation}\label{#1}}
\newcommand{\bea}{\begin{eqnarray}}
\newcommand{\beaa}{\begin{eqnarray*}}
\newcommand{\bean}{\begin{eqnarray}\nonumber}
\newcommand{\beal}[1]{\begin{eqnarray}\label{#1}}
\newcommand{\beadl}[1]{\begin{deqarr}\label{#1}}
\newcommand{\eeadl}[1]{\arrlabel{#1}\end{deqarr}}
\newcommand{\eeal}[1]{\label{#1}\end{eqnarray}}
\newcommand{\eead}[1]{\end{deqarr}}
\newcommand{\eea}{\end{eqnarray}}
\newcommand{\eeaa}{\end{eqnarray*}}

\newcommand{\be}{\begin{equation}}
\newcommand{\ee}{\end{equation}}

\DeclareFontFamily{OT1}{rsfs}{}
\DeclareFontShape{OT1}{rsfs}{m}{n}{ <-7> rsfs5 <7-10> rsfs7 <10->
rsfs10}{} \DeclareMathAlphabet{\mycal}{OT1}{rsfs}{m}{n}

\newcounter{mnotecount}[section]

\newcommand{\N}{{\mathbb N}}

\newcommand{\rmnote}[1]{}%{\mnote{#1}}

\newcommand{\Ric}{\operatorname{Ric}}

\def\mysavedown#1{\edef\mysubs{\mysubs#1}}
\def\mysaveup#1{\edef\mysups{\mysups#1}}
\def\mydown#1{{\mytensor}_{\vphantom{\mysubs}#1}}
\def\myup#1{{\mytensor}^{\vphantom{\mysups}#1}}
\def\tensor#1#2{
  #1
  \def\mytensor{\vphantom{#1}}
  \def\mysubs{\relax}
  \def\mysups{\relax}
  \let\down=\mysavedown
  \let\up=\mysaveup
  #2
  \let\down=\mydown
  \let\up=\myup
  #2
  }

\newcommand{\Hess}{\operatorname{Hess}}
\newcommand{\Ker}{\operatorname{Ker}}
\renewcommand{\Im}{\operatorname{Im}}

\newcommand{\R}{\mathbb R}

\renewcommand{\div}{\operatorname{div}}
\renewcommand{\epsilon}{\varepsilon}
\renewcommand{\hat}{\widehat}

\def\crn#1#2{{\vcenter{\vbox{
        \hbox{\kern#2pt \vrule width.#2pt height#1pt
           }
          \hrule height.#2pt}}}}

\renewcommand{\hbar}{{\overline h}}

\newcommand{\pre}[2]{{{\vphantom{#2}}^{#1}}\kern-.2ex{#2}}

\theoremstyle{plain}
\newtheorem{theorem}{\sc Theorem}[section]

\newtheorem{lemma}[theorem] {\sc Lemma}
\newtheorem{Lemma}[theorem] {\sc Lemma}
\newtheorem{proposition}[theorem]{\sc Proposition}

\newtheorem{corollary}[theorem] {\sc Corollary}
\newtheorem{Corollary}[theorem] {\sc Corollary}

\theoremstyle{definition}
\newtheorem{definition}[theorem]{Definition}

\newtheorem{remark}[theorem]{\sc  Remark\rm}

\numberwithin{equation}{section}

\date{October 8, 2007}

\begin{document}

\title[Spectrum of $\Delta_L$ on A.H. surfaces]
{Spectrum of the Lichnerowicz Laplacian on asymptotically
hyperbolic surfaces}
\author[E.
Delay]{Erwann Delay} \address{Erwann Delay, Laboratoire d'analyse
non lin\'eaire et g\'eom\'etrie, Facult\'e des Sciences, 33 rue
Louis Pasteur, 84000 Avignon, France}
\email{Erwann.Delay@univ-avignon.fr}
\urladdr{http://www.math.univ-avignon.fr/Delay}

\begin{abstract}
We show that, on any asymptotically hyperbolic surface, the
essential spectrum of the Lichnerowicz Laplacian $\Delta_L$ contains
the ray $[\frac{1}{4},+\infty[$. If moreover the scalar curvature is
constant then $-2$ and $0$ are infinite dimensional eigenvalues. If,
 {in addition, the inequality $\langle \Delta u,
u\rangle_{L^2}\geq \frac14||u||^2_{L^2}$ holds for all smooth
compactly supported function $u$,} then there  {is} no other
 {value} in the spectrum.
\end{abstract}

\maketitle

\\
\newline

\tableofcontents
\section{Introduction}\label{section:intro}
This article is a complement of the papers \cite{delay:spectre},
\cite{delay:TT} where the study of the Lichnerowicz Laplacian
$\Delta_L$ is given in dimension $n$ greater than 2. We refer the
reader to those papers for all the motivations. In the preceding
papers, the spectrum was only given for $n\geq3$ because of the
natural relation to the prescribed Ricci curvature problem. In
dimension 2 this study does not appear because the corresponding
problem is conform. The  {present} paper, firstly given for
 {completeness},  {appears} to be particulary interesting
because of  {the quite big differences} with the  {other}
dimensions.

For instance on the hyperbolic space, when $n\geq3$ the spectrum
of $\Delta_L$ on trace free symmetric two tensors is the ray
$$
[\frac{(n-1)(n-9)}{4},+\infty[.
$$
This spectrum is essentially characterized by  non trivial trace free tensors on
the boundary at infinity. In dimension 2 (so 1 at infinity) those
tensors  {do} not  {exist}, and the situation is very
different.

Also, in dimension two, the cohomology of the manifold appears
naturally in the spectrum. This situation was already  {noticed}
by Avez \cite{avez}\cite{avez2} and Buzanca
\cite{buzzanca}\cite{buzzanca2}.

The principal result is the following
\begin{theorem}\label{maintheorem}
Let $(M,g)$  {be an} asymptotically hyperbolic surface. The
essential spectrum of $\Delta_L$ on trace free symmetric two tensors
contains the ray $[1/4,+\infty[$. If moreover $g$ has constant
scalar curvature $R=-2$ then $-2$ and $0$ are also in the essential
 {spectrum. Moreover their eigenspaces are in  {one to one}
correspondance with the space of harmonic one forms respectively in
$L^4$  {and} in $L^2$ (in particular they are infinite
dimensionnal).} Finally, if  {in addition},
 as for the hyperbolic plane, for all smooth
compactly supported function $u$, $\langle \Delta u,
u\rangle_{L^2}\geq \frac14||u||^2_{L^2}$, then  {the} spectrum of
$\Delta_L$ is
$$
\{-2\}\cup\{0\}\cup[\frac{1}{4},+\infty[.$$
\end{theorem}

Along the paper we also obtain some relative results on more
general surfaces, with or without constant scalar curvature.

\noindent{\sc Acknowledgements} I am grateful  {to} N. Yeganefar
for  {discussions} on forms and to F. Gautero for his comments on
the original manuscript.

\section{Definitions, notations and conventions}\label{sec:def}

Let $\overline{M}$ be a smooth, compact surface with boundary
$\partial {M}$. Let $M:=\overline{M}\backslash\partial{M}$  {be}
a non-compact surface without boundary. In our context the boundary
$\partial {M}$ will play the role of a \emph{conformal boundary at
infinity} of $M$. Let $g$ be a Riemannian metric on $M$.  {The
manifold} $(M,g)$ is {\it conformally compact} if there exists on
$\overline{M}$ a smooth defining function $\rho$ for $\partial M$
(that is $\rho\in C^\infty(\overline{M})$, $\rho>0$ on $M$, $\rho=0$
on $\partial {M}$ and $d\rho$  {is} nowhere vanishing on
$\partial M$) such that $\overline{g}:=\rho^{2}g$ is a
$C^{2,\alpha}(\overline {M})\cap C^{\infty}_0(M)$ Riemannian metric
on $\overline{M}$ {. We} will denote by $\hat{g}$ the metric
induced on $\partial M$. Now if $|d\rho|_{\overline{g}}=1$ on
$\partial M$, it is well known (see \cite{Mazzeo:hodge} for
instance) that $g$ has asymptotically sectional curvature $-1$ near
its boundary at infinity {. In this} case we say that $(M,g)$ is
{\it asymptotically hyperbolic}. Along the paper, it will be assumed
sometimes than $(M,g)$ has constant scalar curvature {:} then
 {the} asymptotic hyperbolicity enforces the normalisation
\bel{EE} R(g)=-2\;, \ee where $R(g)$ is the scalar curvature of $g$.

The basic asymptotically hyperbolic surface is the real hyperbolic
Poincar\'e disc. In  {this} case $M$ is the unit disc of $\R^2$,
with the hyperbolic metric
\bel{hyprep} g_0=\omega^{-2}\delta
 \;,
\ee% $$
$\delta$ is the Euclidean metric,
$\omega(x)=\frac{1}{2}(1-|x|_\delta^2)$.

We denote by ${\mathcal T}^q_p$ the set of rank $p$ covariant and
rank $q$ contravariant tensors. When $p=2$ and $q=0$, we denote by
${\mathcal S}_2$ the subset of symmetric tensors, and by
$\mathring{\mathcal S}_2$ the subset of ${\mathcal S}_2$ of trace
free symmetric tensors. We use the summation convention, indices are
lowered with $g_{ij}$ and raised with its inverse $g^{ij}$.

The Laplacian is defined as
$$
\triangle=-tr\nabla^2=\nabla^*\nabla,
$$
where $\nabla^*$ is the $L^2$ formal adjoint of $\nabla$. In
dimension 2, the Lichnerowicz Laplacian acting on trace free
symmetric covariant 2-tensors is
$$
\triangle_L=\triangle+2R,
$$
where $R$ is the scalar curvature of $g$.

For $u$ a covariant 2-tensorfield on $M$ we define the divergence of
$u$ by $$ (\mbox{div}\;u)_i=-\nabla^ju_{ji}.$$ If $u$ is a symmetric
covariant 2-tensorfield on $M$, it can be seen as a one form
with values in the cotangent
bundle. Thus we can define its exterior differential with
$$
(d^\nabla u)_{ijk}:=\nabla_iu_{jk}-\nabla_ju_{ik},
$$
which is a two form with values
the cotangent bundle.

For $\omega$, a one form on $M$, we define its divergence
$$
d^*\omega=-\nabla^i\omega_i,
$$
the symmetric part of its covariant derivative :
$$
({\mathcal
L}\omega)_{ij}=\frac{1}{2}(\nabla_i\omega_j+\nabla_j\omega_i),$$
(note that ${\mathcal L}^*=\mbox{div}$) and the trace free part of
that last tensor :
$$
(\mathring{\mathcal
L}\omega)_{ij}=\frac{1}{2}(\nabla_i\omega_j+\nabla_j\omega_i)+\frac{1}{2}d^*\omega
g_{ij}.$$ The well known \cite{besse:einstein} Weitzenb\"ock
formula for the Hodge-De Rham Laplacian on 1-forms, in dimension
2, reads
$$
\Delta_H\omega_i=\nabla^*\nabla\omega_i+Ric(g)_{ik}\omega^k
=\nabla^*\nabla\omega_i+\frac{R}{2}\omega_i.$$
% All quantities
%relative to $\hat{g}$ will have a hat or will be indexed by
%$\hat{g}$ ($\nabla_{\hat{g}}$, div$_{\hat{g}}$, $d_{\hat{g}}$,
%$\hat{\Delta}_L$,...)
We recall also the Weitzenb\"ock formula
$$
\Delta_K:=(d^\nabla)^*d^\nabla+\div^*\div=\Delta+R=\Delta_L-R.
$$

 For a one form $\omega$, we will consider the trace free symmetric
 covariant two tensor defined by
 $$
 (\mathring{S}\omega)_{ij}=\omega_i\omega_j-\frac{|\omega|^2}{2}g_{ij}.
 $$

A {\it TT-tensor} (Transverse Traceless tensor) is by definition a
symmetric divergence free and trace free  covariant 2-tensor.

$L^2$ denotes the usual Hilbert space of functions or tensors with
the product (resp. norm)
$$
\langle u,v\rangle _{L^2}=\int_M\langle u,v\rangle d\mu_g \mbox{
(resp. }|u|_{L^2}=(\int_M |u|^2d\mu_g)^{\frac{1}{2}}),
$$
where $\langle u,v\rangle $ (resp.  $|u|$) is the  usual product
(resp. norm) of functions or tensors relative to $g$, and the
measure $d\mu_g$ is the usual measure relative to $g$ (we will
omit the term $d\mu_g$). For $k\in\N$, $H^k$ will denote the
Hilbert space of functions or tensors with $k$-covariant
derivative in $L^2$, endowed with its standard  product and norm.

We will first work near the infinity of $M$, so it is convenient
to define for small $\epsilon>0$, the manifold
$$
M_\epsilon=\{x\in M, \rho(x)<\epsilon\}.
$$
It is well know that near infinity, we  {can} choose the defining
function $\rho$ to be the $\overline{g}$-distance to the boundary.
Thus, if $\epsilon$ is small enough, $M_\epsilon$ can be identified
with $(0,\epsilon)\times\partial M$ equipped with the metric
$$
g=\rho^{-2}(d\rho^2+\hat{g}(\rho)d\theta^2),
$$
where $\{\hat{g}(\rho)\}_{\rho\in(0,\epsilon)}$ is a family of
smooth, positive functions on $\partial M$, with
$\hat{g}(0)=\hat{g}$.

Let $P$ be an uniformly degenerate elliptic operator of order 2 on
some tensor bundle  over $M$ (see \cite{GL} for more details). We
recall here a criterion for $P$ to be semi-Fredholm. We first need
the

\begin{definition}\label{defi:estias}
We say that $P$ satisfies the asymptotic estimate
$$
\langle Pu,u\rangle _{L^2} {\geq_\infty }
C||u||^2_{L^2}\;(\mbox{resp. }||Pu||_{L^2} {\geq_\infty } C
||u||_{L^2} )
$$
if for all $\epsilon>0$, there  {exists} $\delta>0$ such that,
for all  {smooth $u$} with compact support in $M_\delta$, we have
$$
\langle Pu,u\rangle _{L^2} {\geq} (C-\epsilon)
||u||^2_{L^2}\;(\mbox{resp. }||Pu||_{L^2} {\geq} (C-\epsilon)
||u||_{L^2}).
$$
\end{definition}

 {Proposition \ref{prop:semifred} below} is standard in the
context of non-compact manifolds (see \cite{delay:spectre} for instance).
It shows that the essential spectrum is characterized near infinity.

\begin{proposition}\label{prop:semifred}
Let $P:H^2\longrightarrow L^2$.  {Then} $P$ is semi-Fredholm (ie.
has finite dimensional kernel and closed range) if and only if $P$
satisfies an asymptotic estimate
$$
||Pu||_{L^2} {\geq_\infty } c ||u||_{L^2}
$$ for some $c>0$.
\end{proposition}

 {This proposition} will be used to compute the essential
spectrum of $\Delta_L$  {which is, by definition,} the closed set
$$
\sigma_{e}(\Delta_L)=\{\lambda\in\Reel,\;\Delta_L-\lambda Id
\mbox{ is not semi-Fredholm} \}.
$$

\section{Commutators of some natural operators}
\begin{lemma}\label{divL} On one forms,
we have $$\div\circ{\mathring{\mathcal
L}}=\frac{1}2(\Delta-\frac{R}{2})=\frac{1}2(\Delta_H-{R}).$$
\end{lemma}
\begin{proof}
In local coordinates,  $2\div\circ{\mathring{\mathcal L}}(\omega)$
is equal to :
\begin{eqnarray*}
-\nabla^{i}(\nabla_i\omega_j+\nabla_j\omega_i-\nabla^k\omega_kg_{ij})&=&
\Delta
\omega_j-\nabla^k\nabla_j\omega_k+\nabla_j\nabla^k\omega_k\\
&=&(\Delta-\Ric)\omega_j\\
&=&(\Delta-\frac{R}{2})\omega_j.
\end{eqnarray*}
\end{proof}
Recall that in dimension 2  {(see Corollary 3.2 of
\cite{delay:TT} for instance) we have:}
\begin{Lemma}
 {Let $(M,g)$ be a Riemannian surface with Levi-Civita connexion
$\nabla$. Then the following equality holds for trace free symmetric
covariants two tensors:}
$$
\div\circ\Delta_L=\Delta_H\circ\div.
$$
\end{Lemma}
So we obtain
\begin{Corollary}
If  {$h$ is a trace free symmetric covariant two tensor} with
$\Delta_L h=\lambda h$ then $\Delta_H \div h=\lambda \div h$.
\end{Corollary}

\begin{Lemma}\label{deltacL}
On a Riemannian surface  with Levi-Civita connexion $\nabla$, we
have
$$
\Delta_L\circ{\mathcal L}={\mathcal
L}\circ\Delta_H-\mathring{S}(dR,.),
$$
where
$\mathring{S}(dR,\xi)_{ij}=\frac{1}{2}(\nabla_jR\xi_{i}+\nabla_iR\xi_{j}-\nabla_pR\xi^pg_{ij})$.
In particular
$$
\Delta_L\circ\mathring{\mathcal L}=\mathring{\mathcal
L}\circ\Delta_H-\mathring{S}(dR,.).
$$
 {Moreover} $R$ is constant iff
$$
 \Delta_L\circ{\mathcal
L}={\mathcal L}\circ\Delta_H \mbox{ so }
\Delta_L\circ\mathring{\mathcal L}=\mathring{\mathcal
L}\circ\Delta_H.
$$
\end{Lemma}
\begin{proof}
The first part comes from \cite{delay:TT} lemma 3.3 where here
$\Ric(g)=\frac{R(g)}{2}g$. Now, if $\Delta_L\circ{\mathcal
L}={\mathcal L}\circ\Delta_H$ then for  {any} one form $\xi$,
$\mathring{S}(dR,\xi)=0$. At any point $x\in M$, we  {take} an
orthonormal basis $(e_1,e_2)$ on $T_x^*M$, and choose
$\xi=e_1$ {. W}e then see that the matrix of
$\mathring{S}(dR,\xi)$  {has} the form
$\left(\begin{array}{cc}a&b\\b&-a\\ \end{array}\right)$, where
$(a,b)$ are the  {coordinates} of $dR$. We finally deduce that
$dR=0$.
\end{proof}
\section{Some decompositions of trace free symmetric two tensors}
In this section, we recall two well known natural decompositions.
We give their simple proofs for completeness.
\begin{lemma}For all $k\in\N$,
$$H^{k+1}(M,\mathring{\mathcal S}_2)=\ker \div\oplus \Im
\mathring{\mathcal L},$$ where the decomposition is orthogonal in
$L^2$.
\end{lemma}
\begin{proof}
For $\omega\in C^\infty_c(M)$, we have
$$
\int_M<\mathring{\mathcal L}(\omega),h>=\int_M<\omega,\div h>.
$$
Thus $\mathring{\mathcal L}^*=\div$,  {which gives} $(\Im
\mathring{\mathcal L})^\bot=\Ker \div$.

\end{proof}

%
%We will be interested by smooth symmetric trace free covariant two
%tensors $h$ in $L^2$ such that theirs  image from some operators
%of the type $\Delta+$ curvatures terms are in $L^2$. From standard
%elliptic regularity, it easy to see that they are in $H^{k+1}$,
%for all integer $k\geq -1$\erw{comme ca c'est faux}. The preceding
%corollary then shows they can be decomposed in a unique way by
%\bel{decompoh} h=\mathring{h}+\mathring{\mathcal L}\omega, \ee
%where $\mathring{h}$ and $\omega$ are in $H^\infty$ (then smooth
%and in $L^2$) and $\div\mathring{h}=0$.

\begin{lemma}\label{decompoform}
For all $k\in\N$,
$$H^{k+1}(M,{\mathcal T}_1)=\ker \Delta_H\oplus \Im d\oplus \Im (*d),$$ where the
decomposition is orthogonal in $L^2$.
\end{lemma}
\begin{remark}
Recall that, from the definition of $\Delta_H$, we have:
$\ker\Delta_H=\ker d\cap\ker d^*.$
\end{remark}

\begin{proof}
First, from the definition of $d^*$, it is clear that $(\Im
d)^\bot=\ker d^*$, and so $H^{k+2}(M,{\mathcal T}_1)=\ker d^*\oplus
\Im d$. For all $H^1$ function $u$ and all $H^1$  {one forms}
$\omega$, we have
$$
\int_M \langle *du,\omega\rangle=\int_M \langle
d^**u,\omega\rangle=\int_M\langle *u,d\omega\rangle=\int_M \langle
u,*d\omega\rangle.
$$
 {As a consequence,} if $\langle *du,\omega\rangle_{L^2}=0$ for
all $u\in C^\infty_c(M)$, then $d\omega=0$, and if in addition
$d^*\omega=0$ then $\Delta_H\omega=0$. This  {shows} that $\ker
d^*=\Im(*d)\oplus\ker \Delta_H$.
\end{proof}

>From  {Lemma} \ref{decompoform}, any one form $\omega$ in $H^1$
can be decomposed in a unique way with \bel{decompoo}
\omega=\eta+du+*dv,\ee where $\Delta_H\eta=0$.

\section{The spectrum on TT-tensors}

\begin{Lemma}\label{TTcarac}
 {Let $M$ be any Riemannian surface. If $h\in
C^2(M,\mathring{\mathcal S}_2)$,
then the following properties are equivalent:}\\
(i)  $\div h=0$,\\
(ii) $d^\nabla h=0$,\\
(iii) $h=\mathring{S}(\omega)$, where $\omega$ is a harmonic one
form.\\
They imply  \\
(iv) $\Delta_L h=R h$.\\
Moreover, if $h\in L^2$, then (iv) implies (i), (ii) and (iii).
\end{Lemma}
\begin{proof}
The first part is due to Avez  (\cite{avez} Lemma A and Lemma C).
The second part is simply due to the following Weitzenb\"ock
formula \cite{Koiso}:
$$
(d^\nabla)^*d^\nabla+\div^*\div=\Delta_K=\Delta_L-R,
$$
and the fact that if $h\in L^2$ solves (iv) weakly, then elliptic
regularity gives $h\in H^\infty\subset C^{\infty}$.
\end{proof}
\begin{corollary}
There  {exists} a non trivial eigen-TT-tensor of $\Delta_L$ iff
$R$ is constant. In  {this} case any TT-tensor is  {an
eigentensor} with eingenvalue $R$.
\end{corollary}
\begin{proof}
The "if" part is clear. For the "only if"  {direction}, assume
that $h$ is a non trivial eigen-TT-tensor of $\Delta_L$,  {so
that $\div h =0$ and $\Delta_Lh=\lambda h$ hold for some
$\lambda\in\R$}. From  {Lemma \ref{TTcarac}, $(\Delta_L-R)h=0$
and then $(R-\lambda)h=0$}. If $R\neq\lambda$  near a point, then
$h$ has to be trivial near this point, so from  {the} unique
continuation property, $h$ is trivial. This  {contradicts} the
assumption on $h$ and  {proves} the result.
\end{proof}

\section{Spectrum on $\Im \mathring{\mathcal L}$}

If $R$ is constant then from  {Lemma \ref{deltacL}} and the fact
that $\Delta_H$  {preserves} the decomposition \ref{decompoo}, it
 {suffices} to study the spectrum of $\Delta_L$ on $\Im
\mathring{\mathcal L}$, restricted successively to $\Ker \Delta_H$,
$\Im d$ and $\Im (*d)$.

\begin{lemma}\label{kerdeltaH} When $R$ is constant then
$\mathring{\mathcal L}(\Ker\Delta_H)$ is in the kernel of
$\Delta_L$ {. If} $R$ is moreover negative, $\mathring{\mathcal
L}(\Ker\Delta_H)$ is in one to one correspondance with
$\Ker\Delta_H$.
\end{lemma}
 \begin{proof}
 If
$h=\mathring{\mathcal L}\eta$, with $\eta\in\Ker \Delta_H$, then
from Lemma \ref{deltacL} $\Delta_L h=0$. Now from Lemma \ref{divL}
we have
$$
2\div\circ{\mathring{\mathcal
L}}(\eta)=(\Delta-\frac{R}{2})\eta.$$ Thus if $R<cte<0$, then
$\mathring{\mathcal L}$ is injective on $H^2$.
\end{proof}

We are now interested in the spectrum on $\Im\mathring{\mathcal
L}\circ d$. We begin with  {a} lemma.

\begin{lemma} \label{normeho} If $h=\mathring{\mathcal L}\omega$, with $\omega\in
H^1$ then:
$$
||h||^2_{L^2}=\frac12(||\omega||^2_{H^1}-\int_M(\frac{R}{2}+1)|\omega|^2).
$$
In particular, if $R=-2$ we obtain $2||h||^2=||\omega||^2_{H^1}$.
\end{lemma}

\begin{proof}
Using Lemma \ref{deltacL} we  {compute:}
\begin{eqnarray*}
\int_M|h|^2&=&\langle\mathring{\mathcal L}\omega,\mathring{\mathcal L}\omega\rangle_{L^2}\\
&=&\langle \div \mathring{\mathcal L}\omega,\omega\rangle_{L^2}\\
&=&\frac12\int_M(|\nabla\omega|^2-\frac{R}{2}|\omega|^2)\\
&=&\frac12(||\omega||^2_{H^1}-\int(\frac{R}{2}+1)|\omega|^2).
\end{eqnarray*}
\end{proof}

\begin{corollary}\label{normeho2} On an A.H. surface, for all $\epsilon>0$, there {exists}
$\delta_0>0$ small such that, for all $\delta\in(0,\delta_0)$ and
all one  {forms} $\omega$ with compact support in $M_\delta$,
 {if} $h=\mathring{\mathcal L}\omega$ then
$$||\omega||^2_{H^1}\geq 2(1-\epsilon)||h||^2_{L^2}.$$
\end{corollary}

\begin{proof}
$$
\begin{array}{lll}
\displaystyle
||\omega||^2_{H^1}-\int_M(\frac{R}{2}+1)|\omega|^2&=&
\displaystyle
 ||\omega||^2_{H^1}-\int_{M_\delta}
O(\rho)|\omega|^2\\&\leq& \displaystyle
||\omega||^2_{H^1}+C\delta||\omega||^2_{L^2}\\&\leq& \displaystyle
\left(1+C\delta\right)||\omega||^2_{H^1},
\end{array}
$$
where $C$ is  a positive constant. Lemma   \ref{normeho}
 {concludes} the proof.

\end{proof}

Let us recall a well known lemma.
\begin{lemma}\label{uo} Let $u$ be a smooth compactly supported function.
If $$\langle\Delta u,u\rangle_{L^2}\geq c||u||^2_{L^2},$$ then
$$\langle\Delta_H du,du\rangle_{L^2}\geq c||du||^2_{L^2}$$ and
$$\langle\Delta_H (*du),(*du)\rangle_{L^2}\geq c||*du||^2_{L^2}.$$
\end{lemma}
\begin{proof}
\begin{eqnarray*}\langle\Delta_H du,du\rangle&&=\langle dd^*du,du\rangle=\langle
d^*du,d^*du\rangle\\&&=||\Delta u||^2\geq c||u||||\Delta u||\geq
c\langle u,\Delta u\rangle=c||du||^2.\end{eqnarray*}
\begin{eqnarray*}\langle\Delta_H (*du),(*du)\rangle&&=\langle d^*d
(*du),*du\rangle=\langle d*du,d*du\rangle\\&&=||\Delta u||^2\geq
c||u||||\Delta u||\geq c\langle u,\Delta
u\rangle=c||du||^2=c||*du||^2.\end{eqnarray*}
\end{proof}
We would like an equivalent to this lemma  {when substituting one
forms to functions. This is achieved} by the following lemma and its
corollary.

\begin{lemma}\label{inegah} Let $\omega$ be a smooth compactly
supported one form. If $$\langle
\Delta_H\omega,\omega\rangle_{L^2}\geq c ||\omega||^2_{L^2}$$ then
$$\langle\Delta_L\mathring{\mathcal L}\omega,\mathring{\mathcal
L}\omega\rangle\geq
\frac{c}2||\omega||^2_{H^1}+\frac12c\int_M(\frac{R}{2}+1)|\omega|^2
-\int_M(\frac R2+1)\langle
\Delta_H\omega,\omega\rangle-\langle\mathring{S}(dR,\omega),\mathring{\mathcal
L}\omega\rangle_{L^2} .$$
\end{lemma}

\begin{proof}
\begin{eqnarray*}
\langle\Delta_L \mathring{\mathcal L}\omega, \mathring{\mathcal
L}\omega\rangle_{L^2}&=&\langle\mathring{\mathcal L}\Delta_H
\omega, \mathring{\mathcal
L}\omega\rangle_{L^2}-\langle\mathring{S}(dR,\omega),\mathring{\mathcal
L}\omega\rangle_{L^2} \\
&=&\langle\Delta_H \omega, \div\mathring{\mathcal
L}\omega\rangle_{L^2}-\langle\mathring{S}(dR,\omega),\mathring{\mathcal
L}\omega\rangle_{L^2}
 \\
&=&\frac12\langle\Delta_H \omega,
(\Delta_H-R)\omega\rangle_{L^2}-\langle\mathring{S}(dR,\omega),\mathring{\mathcal
L}\omega\rangle_{L^2}
 \\
&=&\frac12||\Delta_H \omega||_{L^2}^2-\frac{1}2\langle R\Delta_H
\omega,\omega\rangle_{L^2}-\langle\mathring{S}(dR,\omega),\mathring{\mathcal
L}\omega\rangle_{L^2} \\
&\geq&\frac12c||\Delta_H
\omega||_{L^2}||\omega||_{L^2}-\frac{1}2\langle R\Delta_H
\omega,\omega\rangle_{L^2}-\langle\mathring{S}(dR,\omega),\mathring{\mathcal
L}\omega\rangle_{L^2} \\
&\geq&\frac12c\langle\Delta_H \omega,\omega\rangle_{L^2}
-\frac{1}2\langle R\Delta_H
\omega,\omega\rangle_{L^2}-\langle\mathring{S}(dR,\omega),\mathring{\mathcal
L}\omega\rangle_{L^2}\\
&\geq&\frac12c||\nabla\omega||_{L^2}^2+\frac12c\int_M\frac{R}{2}|\omega|^2
-\frac{1}2\langle
R\Delta_H
\omega,\omega\rangle_{L^2}-\langle\mathring{S}(dR,\omega),\mathring{\mathcal
L}\omega\rangle_{L^2}\\
&\geq&\frac12c||\nabla\omega||_{L^2}^2+\frac12c\int_M\frac{R}{2}|\omega|^2
+\langle \Delta_H \omega,\omega\rangle_{L^2} \\& &\;\;\
-\frac{1}2\langle (R+2)\Delta_H
\omega,\omega\rangle_{L^2}-\langle\mathring{S}(dR,\omega),
\mathring{\mathcal L}\omega\rangle_{L^2}\\
&\geq&\frac12c||\nabla\omega||_{L^2}^2+\frac12c\int_M\frac{R}{2}|\omega|^2
+c||\omega||^2_{L^2} \\& &\;\;\ -\frac{1}2\langle (R+2)\Delta_H
\omega,\omega\rangle_{L^2}-\langle\mathring{S}(dR,\omega),
\mathring{\mathcal L}\omega\rangle_{L^2}\\
&\geq&\frac12c||\omega||_{H^1}^2+\frac12c\int_M(\frac{R}{2}+1)|\omega|^2
\\& &\;\;\ -\frac{1}2\langle (R+2)\Delta_H
\omega,\omega\rangle_{L^2}-\langle\mathring{S}(dR,\omega),
\mathring{\mathcal L}\omega\rangle_{L^2}.\\
\end{eqnarray*}

\end{proof}

\begin{remark}
 {Under the assumptions of Lemma \ref{uo}, the assumptions} of
lemma \ref{inegah} are satisfied by $\omega=du$ or $\omega=*du$.
\end{remark}

Proposition \ref{inegah} together with Lemma \ref{normeho} give:

\begin{corollary}\label{estideltaL}
If $R=-2$ and $\langle \Delta_H\omega,\omega\rangle_{L^2}\geq c
||\omega||^2_{L^2}$ then $$\langle\Delta_L\mathring{\mathcal
L}\omega,\mathring{\mathcal L}\omega\rangle_{L^2}\geq
c||\mathring{\mathcal L}\omega||^2_{L^2}.$$
\end{corollary}

In the A.H.  {setting} we have
\begin{corollary}\label{estideltaLinfty} On an A.H. surface, for $\omega\in \Im
d$ or $\omega\in\Im *d$, $$\langle\Delta_L\mathring{\mathcal
L}\omega,\mathring{\mathcal L}\omega\rangle_{L^2}\geq_\infty
\frac14||\mathring{\mathcal L}\omega||^2_{L^2}.$$
\end{corollary}

\begin{proof}

It is well know that on A.H. surfaces,  {$\langle\Delta
u,u\rangle\geq_{\infty}\frac14||u||^2_{L^2}$ holds.} Then (see lemma
\ref{uo} for instance)
$\langle\Delta_H\omega,\omega\rangle_{L^2}\geq_\infty\frac{1}{4}|\omega|^2_{L^2}.$
We will show that the  {three} terms in the  {right-hand} side
of Lemma \ref{inegah} do not contribute at infinity. We work with a
one form $\omega$ compactly supported in $M_\delta$ with small
$\delta$. We recall that   $R+2=O(\rho)$ and $||d(R+2)||=O(\rho).$

Let us begin with \begin{eqnarray*} \int_M(\frac R2+1)\langle
\Delta_H\omega,\omega\rangle&=&\int_M(\frac R2+1)\langle
\Delta\omega,\omega\rangle+\int_M(\frac
R2+1)\frac{R}2||\omega||^2\\&=& \int_M\omega^j\nabla^i(\frac R2+1)
\nabla_i\omega_j+\int_M(\frac R2+1)
||\nabla\omega||^2\\&&\;\;+\int_M(\frac
R2+1)\frac{R}2||\omega||^2.\\
\end{eqnarray*}
The two last terms are clearly bounded in absolute value by
$C_1\delta||\omega||^2_{H_1}$. Let
$A(\omega)=\int_M\omega^j\nabla^i(R+2) \nabla_i\omega_j$.
 {Then:}
\begin{eqnarray*}
|A(\omega)|&\leq& ||\nabla\omega||_{L^2}||
d(R+2)\omega||_{L^2}\\&\leq&||\nabla\omega||_{L^2}(\int_M||d(R+2)||^2||\omega||^2)^{1/2}\\
&\leq&C_2\delta||\nabla\omega||_{L^2}(\int_M||\omega||^2)^{1/2}\\
&\leq&\frac {C_2}2\delta||\omega||^2_{H^1}.\\
\end{eqnarray*}
 {We so get:}
$$
|\int_M(\frac R2+1)\langle \Delta_H\omega,\omega\rangle|\leq
C_3\delta||\omega||^2_{H^1}.
$$
 {The term $\langle\mathring{S}(dR,\omega),\mathring{\mathcal
L}\omega\rangle_{L^2}$ proceed in a manner similar to $A(\omega)$ to
obtain the same   estimate, perhaps with a different constant}. Finally, the term
$\int_M(\frac{R}{2}+1)|\omega|^2$ is clearly bounded in absolute
value by $C_4\delta||\omega||^2_{H_1}$. The conclusion  {follows}
from Lemma \ref{inegah}, the  {triangular} inequality and
Corollary \ref{normeho2}.

\end{proof}

\begin{proposition}\label{prop:noestias}
 {Let $\lambda\geq\frac{1}{4}$ and $C>0$. Let $P$ be the operator
$\Delta_L- \lambda Id : H^2 \longrightarrow L^2$. There is no
asymptotic estimate
$$
|Pu|_{L^2}{\geq_\infty } C |u|_{L^2},
$$
for $P$ on $\Im (\mathring{\mathcal L}\circ d)$.}
\end{proposition}
\begin{proof}
Let $\lambda\geq\frac{1}{4}$, and
 $\mu:=\sqrt{\lambda-\frac{1}{4}}$.
The idea of the proof is to construct a family of tensors
$\{h_R\}=\{\mathring{\mathcal L}(df_R)\}=\{\mathring{\Hess}f_R\}$
with compact support in $M_{e^{-R/2}}$ such that $|Ph_R|_{L^2(M)}$
goes to zero when
 $R$ goes to infinity but  $|h_R|_{L^2(M)}$ goes to infinity when $R$ goes
 to infinity.

It is well known (see [L1, lemma 5.1] for example) that we can
change the defining function $\rho$ into a defining function $r$
such that the metric takes the form
$$
g=r^{-2}\overline{g}=r^{-2}(dr^2+\hat{g}(r)),
$$
on $M_{\delta}=]0,\delta[\times \partial_\infty M$ (reducing
$\delta$ if necessary), where $\hat{g}(r)$ is a metric on
$\{r\}\times \partial_\infty M$.

 The non trivial Christoffel symbols of $g=r^{-2}[dr^2+\hat{g}(r)d\theta^2]$
are
$$
\Gamma^r_{rr}=-r^{-1},
$$
$$
\Gamma^r_{\theta\theta}=-\frac{r^2}{2}(-2r^{-3}\hat{g}
+r^{-2}\hat{g}')=r^{-1}\hat{g}-\frac12\hat{g}',
$$
$$
\Gamma^\theta_{\theta r}=\frac{1}{2}(-2r^{-1}
+\hat{g}^{-1}\hat{g}')=-r^{-1}+\frac12\hat{g}^{-1}\hat{g}',
$$
where the primes denote $r$-derivatives. If $f$ is a "radial"
function, ie $f=f(r)$, we  {compute} :
$$
\Hess
f=(f''+r^{-1}f')dr^2+(-r^{-1}\hat{g}+\frac12\hat{g}')f'd\theta^2.
$$
 {We deduce}
$$
\Delta f=-r^2(f''+\frac12\hat{g}^{-1}\hat{g}'f').
$$
We also have
$$
\mathring{\Hess}
f=(\frac12f''+r^{-1}f'-\frac14\hat{g}^{-1}\hat{g}'f')(dr^2-\hat{g}d\theta^2)
=:F_f(r)(dr^2-\hat{g}d\theta^2).
$$
This tensor is in the set ${\mathcal V}_2$ of \cite{GL} page 201:
 {substitute $f$ there by $F_f$ here, $\overline{q}$ there by
$dr^2-\hat{g}d\theta^2$ here, $\rho$ there by $r$ here and the
dimension $n+1$ there by $2$ here}. Recall that the Lichnerowicz
Laplacian in our context is $\Delta_L=\Delta+{\mathcal K}$, where
${\mathcal K}=-4+O(r)$. Thus, from \cite{GL} Lemma 2.9 page 202, we
obtain
$$
(\Delta_L-\lambda)(F(r)\overline{q})=I_2(F(r))\overline{q}+rX(F),
$$
where $$I_2(F)=-r^2F''-4rF'-2fF$$ has for characteristic exponents
$$
s_1,s_2=\frac12(-3\pm\sqrt{1-4\lambda}),
$$
and $X=\overline{a}r^2\frac{d^2}{dr^2}
+\overline{b}r\frac{d}{dr}+\overline{c}$ is a second order
operator polynomial in $r\frac{d}{dr}$ with $\overline{g}$-bounded
coefficients depending on $\overline{g}$ and $\overline{q}$.

In particular, if $\lambda\geq\frac14$ and
$f(r)=\sqrt{r}(a\cos(\mu\ln (r))+b\sin(\mu\ln (r))$, where
$\mu=\sqrt{\lambda-\frac14}$, then $F_f(r)={r}^{-3/2}(A\cos(\mu\ln
(r))+B\sin(\mu\ln (r))+O(r^{-1/2})$ and $(A,B)\neq 0$ if $(a,b)\neq
0$. Thus we obtain
$$
I_2(F_f(r))=O(r^{-1/2}).
$$

 {Let us now define} the function
$$
f_R(r)=f(r)\Psi_R(r),
$$
where $\Psi_R$ is  {as in} Lemma \ref{lem:cutoff}. A simple
calculation  {shows} that
$$
F_{f_R}(r)=\Psi_R(r) F_{f}(r)+O(R^{-1})O(r^{-3/2}).
$$
Therefore:
$$
I_2(F_{f_R}(r))=O(r^{-1/2})+O(R^{-1})O(r^{-3/2}),
$$
and
$$
rX(F_{f_R}(r))=O(r^{-1/2}).
$$
Then:
$$
(\Delta_L-\lambda)(F_{f_R}(r)\overline{q})=O(r^{-1/2})+O(R^{-1})O(r^{-3/2}).
$$

We deduce that
$$
||(\Delta_L-\lambda)(F_{f_R}(r)\overline{q})||^2_{L^2}=O(R^{-1}).
$$
 {On} the other hand, we have
$$
||(F_{f_R}(r)\overline{q})||^2_{L^2}\geq c R,
$$
where $c$ is a positive constant. Letting $R$ going to infinity,
this  {concludes} the proof of the proposition.

\end{proof}

\section{Conclusion}
>From Proposition \ref{prop:noestias} and Corollary
\ref{estideltaLinfty}, the essential spectrum of $\Delta_L$
restricted to $\Im (\mathring{\mathcal L}\circ d)$ is
$$
[\frac{1}{4},+\infty[.
$$
In particular this ray is in the essential spectrum of $\Delta_L$.

If $R$ is constant then the A. H. condition forces $R=-2$. Lemma
\ref{kerdeltaH} shows that any tensor in $\mathring{\mathcal
 L}(\ker\Delta_H)$ is in the kernel of $\Delta_L$.
 The eigenspace for $0$ is then  infinite dimensionnal as $\ker \Delta_H$
(recall that $\mathring{\mathcal L}$ is injective if $R<0$).

 {From Lemma \ref{TTcarac}, any TT-tensor $h$ is an eigentensor
for the eigenvalue $-2$ and there is a harmonic one form omega such
that $h=\mathring{\mathcal S}(\omega)$. Moreover $\omega$ is in
$L^4$ iff $h$ is in $L^2$.}

 {Assume now that
$$\langle\Delta u,u\rangle_{L^2}\geq \frac14 ||u||_{L^2}$$
holds for all smooth compactly supported functions $u$. Then Lemma
\ref{uo} and Corollary \ref{estideltaL} give, when $\omega=df$ or
$\omega=*df$,
$$\langle\Delta_L \mathring{\mathcal L}\omega,\mathring{\mathcal L}\omega\rangle_{L^2}\geq
 \frac14 ||\mathring{\mathcal L}\omega||_{L^2}.$$
 This proves that there are no eigentensors with eigenvalue less than $\frac14$ in $\Im(\mathring{\mathcal
 L}\circ d)$ nor in $\Im(\mathring{\mathcal
 L}\circ (*d))$.
Recall that, when $R$ is a constant, the Lichnerowicz Laplacian
commute with the Hodge Laplacian, and also that the Hodge Laplacian
preserves the decomposition \ref{decompoo}. We so get that on}
$\Im\mathring{\mathcal
 L}$, the essential spectrum of $\Delta_L$ is
 $$
 \{0\}\cup[\frac14,+\infty[.
 $$

This concludes the proof of the main theorem \ref{maintheorem}.

\section{Appendix : a family of cutoff functions}

In this appendix, we give a family of cutoff functions. Standard in
the A.H. context, they can be found in \cite{andersson:elliptic},
Definition 2.1 p.1362 for instance.
\begin{lemma}\label{lem:cutoff}
Let $(M,g,\rho)$ be an asymptotically hyperbolic manifold. For
 $R\in\Reel$ large enough, there exits a cutoff function
$\Psi_R:M\rightarrow[0,1]$ depending only on $\rho$,  supported in
the annulus $\{e^{-8R}<\rho<e^{-R}\}$, equal to $1$ in
$\{e^{-4R}<\rho<e^{-2R}\}$  and which satisfies for $R$ large :
$$
|\frac{d^k\Psi_R}{d\rho^k}(\rho)|\leq \frac{C_k}{R\rho^k},
$$
for all $k\in \Nat\backslash\{0\}$, where $C_k$ is independent of
$R$.
\end{lemma}
\begin{proof}
Let $\chi:\Reel\longrightarrow [0,1] $ be a smooth function equal
to $1$ on $]-\infty,1]$ and  $0$ on $[2,+\infty[$. We define
$$
\chi_R(x):=\chi(\frac{\ln(\rho(x))}{-R}),
$$
we then have $\chi_R:M\longrightarrow [0,1]$  is equal to $1$ on
$\rho\geq e^{-R}$ and $0$ on $\rho\leq e^{-2R}$. Now we define
$$
\Psi_R:=\chi_{4R}(1-\chi_R)
$$
which satisfies the   {announced} properties.
\end{proof}

\bibliographystyle{amsplain}

\bibliography{../references/erwbiblio,%
../references/newbiblio,%
../references/reffile,%
../references/bibl,%
../references/hip_bib,%
../references/newbib,%
../references/PDE,%
../references/netbiblio,%
stationary}
\end{document}